\documentclass[12pt]{article}
\usepackage{amsmath} 
\usepackage[english]{babel}
\usepackage{amsbsy}
\usepackage{amsfonts}
\usepackage{amssymb}
\usepackage{amsthm}
\usepackage[pdftex]{graphicx}
\usepackage[normalem]{ulem}
\usepackage{color}
\usepackage{float}
\usepackage{framed}
\usepackage[hyphens]{url}
\usepackage{mathtools}
\usepackage{color, colortbl}
\usepackage{xcolor}
\usepackage{mathrsfs}
\usepackage{geometry}
\usepackage[T1]{fontenc}
\usepackage[utf8]{inputenc}
\usepackage{authblk}
\usepackage{cases}
\usepackage[numbers]{natbib}
\usepackage{makecell}
\usepackage{algorithmic}
\usepackage{algorithm}
\usepackage{caption}
\usepackage{longtable}
\usepackage{siunitx}
\usepackage{subcaption}
\usepackage{multirow}
\usepackage{booktabs,tabularx}

\usepackage{hyperref}
\hypersetup{
	colorlinks=true,
	linkcolor=blue,
	filecolor=blue,
	citecolor = blue,      
	urlcolor=cyan,
}

\title{Solving the Nonlinear Westervelt Equation in Ultrasound Propagation Using Physics-Informed Neural Networks}

\author{Luke Oluwaseye Joel \thanks{Corresponding author: oluwaseyejoel@gmail.com; ljoel@uj.ac.za}}
\affil[1]{Department of Mathematics and Applied Mathematics\\
	University of Johannesburg\\
	South Africa}\date{}

\begin{document}
	
	\maketitle

	\begin{abstract}
    Modeling the nonlinear acoustic wave propagation in biological tissues needs accuracy, especially in ultrasound imaging and therapeutic applications. The 1D Westervelt equation, $p_{xx} - p_{tt} + \gamma
    p_{ttt} + \eta(p^2)_{tt}=0$, could be used to describe a finite‑amplitude ultrasound propagation in these tissues. This equation accounts for nonlinearity, $\eta$ and attenuation, $\gamma$. To solve this equation, we observe that traditional solvers may struggle to do so correctly. This is due to the difficulties arising from the nonlinear nature of the equation. We solve the equation by defining three synthetic cases - linear ( $\eta=0.001$, $\gamma=10^{-6}$), shock-dominated ( $\eta=0.02$, $\gamma=10^{-6}$), and dissipation-dominated ( $\eta=0.01$, $\gamma=5\times10^{-4}$), and using Physics‑Informed Neural Networks (PINNs) to showcase the solution. A sinusoidal activation function is used to capture the oscillatory behavior of the problem. The results of our study, with PDE residual converging to at least $10^{-6}$, bundary condition loss to $10^{-5}$ and the PINN loss converges to $10^{-5}$ for each case at $30 \, 000$ epochs, show that $\gamma$  scales the rate of amplitude decay), while $\eta$ governs frequency-dependent phase speeds.
	\end{abstract}
	
	
	\section{Introduction}
	Nonlinear acoustic wave propagation \cite{westervelt1963rotational, nikolic2019priori} is central to many modern scientific and engineering fields. These include high‑intensity focused ultrasound (HIFU) \cite{ter2007high,elhelf2018high,evans2007high,filonenko2001effect,pinton2011effects}, medical imaging \cite{ward1997nonlinear,shishlenin2024nonlinear,sarvazyan2013acoustic} , lithotripsy \cite{howle1998lithotripsy, coleman1993review,cleveland2005modeling}, nondestructive testing \cite{donskoy2001nonlinear,nasrollahi2017nondestructive}, and underwater acoustics \cite{yuldashev2025numerical,campo2024signal}. A mathematical model for finite‑amplitude sound propagation in thermoviscous fluids is described the Westervelt equation \cite{westervelt1963rotational, shevchenko2015absorbing,shevchenko2015ab},
	
	\begin{equation}
		\frac{\partial^2 p}{\partial t^2} - c_0^2 \nabla^2 p - b \nabla^2 \frac{\partial p}{\partial t} = \frac{\beta}{\rho_0 c_0^2} \frac{\partial^2 p^2}{\partial t^2},
		\label{eq:westervelt_general}
	\end{equation}
	where $p$ is the pressure, $c_0$ is the small‑signal sound speed, $b$ is thermoviscous diffusivity, and $\beta$ is the nonlinearity coefficient. For one‑dimensional motion along the $x$-axis, the equation is given by
	
	\begin{equation}
		\frac{\partial^2 p}{\partial x^2} - \frac{1}{c_0^2}\frac{\partial^2 p}{\partial t^2} + \frac{\delta}{c_0^4}\frac{\partial^3 p}{\partial t^3} + \frac{\beta}{\rho_0 c_0^4}\frac{\partial^2 (p^2)}{\partial t^2} = 0,
		\label{eq:westervelt_1d}
	\end{equation}

The third term $\frac{\delta}{c_0^4}\frac{\partial^3 p}{\partial t^3}$ is the frequency‑dependent attenuation while the fourth  term $\frac{\beta}{\rho_0 c_0^4}\frac{\partial^2 (p^2)}{\partial t^2}$ is the quadratic nonlinearity. This latter term is responsible for harmonic generation and waveform distortion. The coefficient $\beta$ relates to the standard nonlinearity parameter $B/A$ through $\beta = 1 + B/(2A)$.

The initial conditions are
	\begin{equation}
		p(x, 0) = 0, \qquad p_t(x, 0) = 0, \qquad p_{tt}(x, 0) = 0, \quad \forall x \in [0, L],
		\label{eq:ic}
	\end{equation}
A pulsed ultrasound source is applied at the transducer face. To avoid numerical artifacts from sharp discontinuities, we use a smooth Gaussian‑enveloped cosine pulse. This models a physically realistic transducer response
	\begin{equation}
		p(0, t) = A_p \exp\left(-\frac{(t - t_c)^2}{2\sigma^2}\right) \cos\left[\omega_0 (t - t_c)\right].
		\label{eq:bc_source}
	\end{equation}
Here, $A_p$ is pressure amplitude, $\omega_0$ is the fundamental angular frequency, $t_c$ is pulse center time, and $\sigma$ is pulse width (the standard deviation of the Gaussian envelope). This source model is common in nonlinear acoustics. It represents transducers with Gaussian beam profiles and agrees well with experimental measurements of pulsed ultrasound fields.
	
At the far end of the domain ($x = L$), we apply an absorbing (radiation) boundary condition. This minimizes unphysical reflections
	
	\begin{equation}
		p_t(L, t) + c_0 \frac{\partial p}{\partial x}(L, t) = 0.
		\label{eq:bc_absorbing}
	\end{equation}
This Sommerfeld‑type condition lets waves exit the domain with minimal reflection. It effectively simulates an infinite or semi‑infinite medium. We now introduce the following parameters	
	\begin{equation}
		\alpha = \frac{1}{c_0^2}, \qquad \gamma = \frac{\delta}{c_0^4}, \qquad \eta = \frac{\beta}{\rho_0 c_0^4}.
		\label{eq:params_definition}
	\end{equation}
Substituting these into equation \eqref{eq:westervelt_1d} and $c_0 =1$ gives this normalized form	
	\begin{equation}
		p_{xx} -  p_{tt} + \gamma p_{ttt} + \eta (p^2)_{tt} = 0.
		\label{eq:normalized_westervelt}
	\end{equation}
Hence, the complete initial‑boundary value problem for the nonlinear Westervelt equation is given as follows.
	\begin{align}
		p_{xx} - p_{tt} + \gamma p_{ttt} + \eta (p^2)_{tt} &= 0, \qquad (x,t) \in [0, L] \times (0, T], \nonumber \\
		p(x, 0) &= 0, \qquad x \in [0, L], \nonumber \\
		p_t(x, 0) &= 0, \qquad x \in [0, L], \nonumber \\
		p_{tt}(x, 0) &= 0, \qquad x \in [0, L], \label{eq:bc_full_absorbing} \\
		p(0, t) &= A_p \exp\left(-\frac{(t - t_c)^2}{2\sigma^2}\right) \cos\left[\omega_0 (t - t_c)\right], \qquad t \in (0, T], \nonumber \\
		p_t(L, t) + c_0 \frac{\partial p}{\partial x}(L, t) &= 0, \qquad t \in (0, T]. \nonumber
	\end{align}
	where the parameters $\gamma$ and $\eta$ are relative to the cases we are studying.
	
PINNs have been applied to ultrasound wave propagation. However, existing studies seems not to be physically complete. Alkhadhr and Almekkawy \cite{alkhadhr2022physics,alkhadhr2023modeling,alkhadhr2023wave} developed a PINNs frameworks for multi-element transducers. The linear wave equation was used for this without attenuation, $\gamma$ or nonlinearity, $\eta$. Soft-soft constraint combinations were found to outperform hard constraint. Li and Pinton \cite{pinton2011effects} developed PINNs to solve a nonlinear wave equation, achieving about $ 12 \times$ speedup. However, the tissue attenuation, $\gamma$, was neglected. Wang et al. \cite{wang2023physics} also solved the transcranial propagation with a linear 2D wave equation with PINNs. Tube acoustic fields were reconstructed by Luan et al. \cite{luan2025acoustic} as a 1D inverse problem. This problem is fundamentally different from the nonlinear attenuating forward problem. Other related works that use PINNs to non-equivalent equations include the Helmholtz equation \cite{schmid2023physics}, guided waves in solids \cite{stone2023physics}, laser ultrasonics \cite{wang2026sparse}, and bimaterial elastodynamics \cite{chibire2026physics}. From these studies, we observe that no study has solved the full Westervelt equation stated above. Also, no direct PINNs counterpart exists that showcase the different cases as we shall do in this study. Hence, the gap this work seeks to fill.

\section{Physics-Informed Neural Network}
\label{pinnss}
We approximate the solution of the Westervelt equation in \eqref{eq:bc_full_absorbing} using PINNs \cite{raissi2018deep, raissi2019physics, klapa2024machine, joel2025solving}. The idea of PINNs is to embed the PDE, along with its initial and boundary conditions directly into the neural network's loss function. Then, we train the network to yield the PDE solution. This is achieved by minimizing both the PDE residual and the violations of the imposed constraints. We use the following network configuration of 3 hidden layers, 64 neurons/layer, and $\sin$ activation function.

A neural network approximation for the acoustic pressure field, $p(x,t)$, is 
\begin{equation}
	p(x,t) \approx p_\theta(x,t),
\end{equation}
where $\theta$ represents the trainable parameters. The network receives spatial and temporal coordinates as input and produces the predicted pressure. Automatic differentiation computes all derivatives required by the governing equation. This allows exact differentiation of the network with respect to its inputs. To enforce the governing equation within this framework, we define the PDE residual from \eqref{eq:normalized_westervelt} as
\begin{equation}
	\mathcal{R}(x,t) =
	p_{xx} - \alpha p_{tt} + \gamma p_{ttt} + \eta (p^2)_{tt}.
	\label{eq:pde_residual}
\end{equation}

The PINNs' objective is to minimize this residual along with any deviations from the initial and boundary conditions. Training proceeds by minimizing a composite loss function. This total loss, Equation \eqref{eqloss}, has four components namely the PDE residual, the initial condition, the source boundary condition, and the absorbing boundary condition losses. We do not assign weight to the losses.
\begin{equation}
	\mathcal{L} =
	\mathcal{L}_{\text{PDE}}
	+
	\mathcal{L}_{\text{IC}}
	+
	\mathcal{L}_{\text{BC,source}}
	+
	\mathcal{L}_{\text{BC,absorbing}}.
	\label{eqloss}
\end{equation}
With each term defined as
\begin{align}
	\mathcal{L}_{\text{PDE}} &= \frac{1}{N_f} \sum_{i=1}^{N_f} \left|\mathcal{R}(x_i, t_i)\right|^2, \label{eq:loss_pde} \\
	\mathcal{L}_{\text{IC}} &= \frac{1}{N_i} \sum_{i=1}^{N_i} \left( \left|p_{\theta}(x_i, 0)\right|^2 + \left|p_{t,\theta}(x_i, 0)\right|^2 + \left|p_{tt,\theta}(x_i, 0)\right|^2 \right), \label{eq:loss_ic} \\
	\mathcal{L}_{\text{BC,source}} &= \frac{1}{N_{b1}} \sum_{i=1}^{N_{b1}} \left|p_{\theta}(0, t_i) - g(t_i)\right|^2, \label{eq:loss_bc_source} \\
	\mathcal{L}_{\text{BC,absorbing}} &= \frac{1}{N_{b2}} \sum_{i=1}^{N_{b2}} \left|p_{t,\theta}(L, t_i) + c_0 \frac{\partial p_\theta}{\partial x}(L, t_i)\right|^2. 
	\label{eq:loss_bc_absorbing}
\end{align}
Where $N_f$, $N_i$, $N_{b1}$, and $N_{b2}$ are the numbers of collocation points. We use $N_f = N_i = N_{b1} = N_{b2} = 600$ in our simulation. The function $g(t_i)$ is the prescribed transducer source
\begin{equation}
	g(t) = A_p \exp\left(-\frac{(t - t_c)^2}{2\sigma^2}\right) \cos(\omega_0 (t - t_c)).
	\label{eq:g_function}
\end{equation}
We optimize the network parameters $\theta$ by minimizing the total loss. This is done with Adam gradient‑based algorithm (learning rate=$10^{-3}$ and $30\, 000$ epochs). During training, collocation points are sampled across the spatio‑temporal domain. This ensures that the PDE residual is minimized globally. Once training converges, the trained PINNs gives a continuous approximation of the acoustic pressure field over the entire domain. It is should noted that this requires no predefined numerical mesh. We implement the neural network using the PyTorch deep learning framework.

\section{Results and Discussion}
\label{result}
This section presents and discusses our results. We solved Equation \eqref{eq:bc_full_absorbing} for three distinct cases. The PINN framework successfully captured the characteristic wave phenomena unique to each case. It also revealed clinically differences in harmonic generation and acoustic energy distribution with depth.

\begin{table}[h]
	\centering
	\caption{Normalized Parameters for the Simulation}
	\label{tab:parameters}
	\begin{tabular}{lcc}
		\toprule
		\textbf{Parameter} & \textbf{Symbol} & \textbf{Value}  \\
		\midrule
		Wave propagation coefficient & $\alpha$ & 1.0  \\
		Sound speed & $c_0$ & 1.0  \\
		Amplitude & $A_p$ & 1.0  \\
		Fundamental frequency & $\omega_0$ & 4.0 \\
		Pulse center time & $t_c$ & 2.0  \\
		Pulse width  & $\sigma$ & $0.8/\sqrt{2}$ \\
		Domain depth & $L$ & 6.0  \\
		Simulation time & $T$ & 8.0  \\
		\bottomrule
	\end{tabular}
\end{table}

\noindent Table \ref{tab:parameters} summarizes the parameters used in this study. These values are normalized and dimensionless for the computational model. The normalization follows standard practice in nonlinear acoustics. Spatial coordinates are scaled by the wavelength, and time is scaled by the period.
\begin{table}[h]
	\centering
	\caption{Parameter Cases for the Simulation}
	\label{tab:cases}
	\begin{tabular}{lcc}
		\toprule
		\textbf{Case}  & \textbf{$\gamma$} & \textbf{$\eta$} \\
		\midrule
		1 &  $1.0\times10^{-6}$ & 0.001 \\
		2 &  $1.0\times10^{-6}$ & 0.02 \\
		3 &  $6.0\times10^{-4}$ & 0.01 \\
		\bottomrule
	\end{tabular}
\end{table}

\noindent Table \ref{tab:cases} lists the cases \cite{purrington2012numerical,connor2002bio,gupta2018numerical,hallaj2001simulations} we use in our study. They span the range of tissue properties found in medical ultrasound. Case~$1$ represents tissues with low attenuation and weak nonlinearity. Examples include water‑like media or homogeneous soft tissues at low acoustic pressures. Case~$2$ represents typical soft tissue properties. Case~$3$ represents tissues with significant attenuation and strong nonlinearity. These include highly attenuating or dense tissues, which are often encountered in high‑intensity therapeutic ultrasound applications.

\subsection{Case 1: Weak Nonlinearity, Low Attenuation}
\begin{figure}[h!]
	\centering
	\begin{subfigure}[t]{0.45\textwidth}
		\centering
		\includegraphics[height=1.65in]{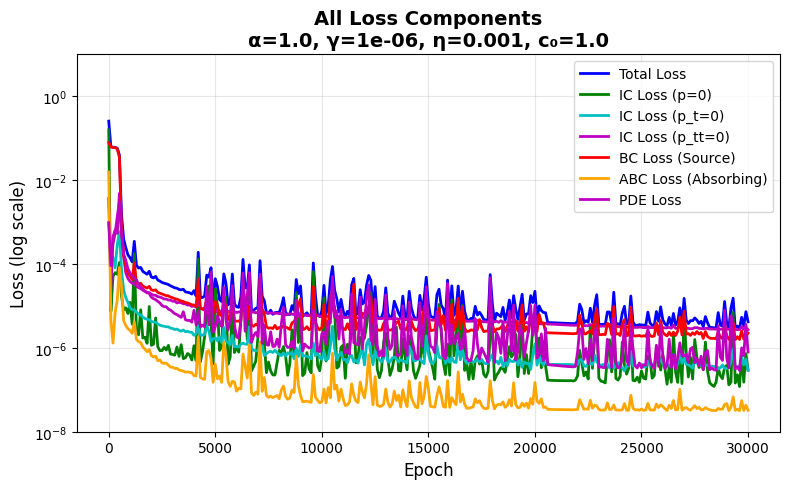}
		\caption{Loss Components}
		\label{fig:case1_loss}
	\end{subfigure}%
	\begin{subfigure}[t]{0.45\textwidth}
		\centering
		\includegraphics[height=1.65in]{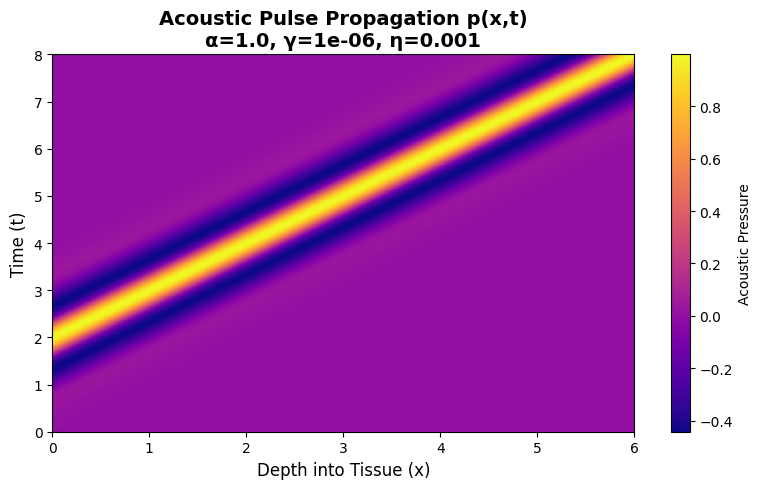}
		\caption{Acoustic Propagation}
		\label{fig:case1_propagation}
	\end{subfigure}
	\caption{Case 1: Loss components and acoustic propagation field.}
	\label{fig:case1}
\end{figure}

\begin{figure}[h!]
	\centering
	\begin{subfigure}[t]{0.45\textwidth}
		\centering
		\includegraphics[height=1.65in]{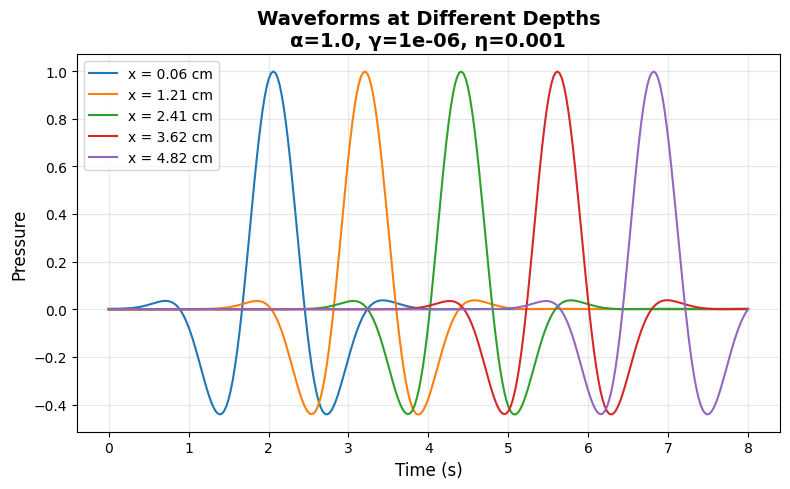}
		\caption{Waveform at Different Depths}
		\label{fig:case1_depth}
	\end{subfigure}%
	\begin{subfigure}[t]{0.45\textwidth}
		\centering
		\includegraphics[height=1.65in]{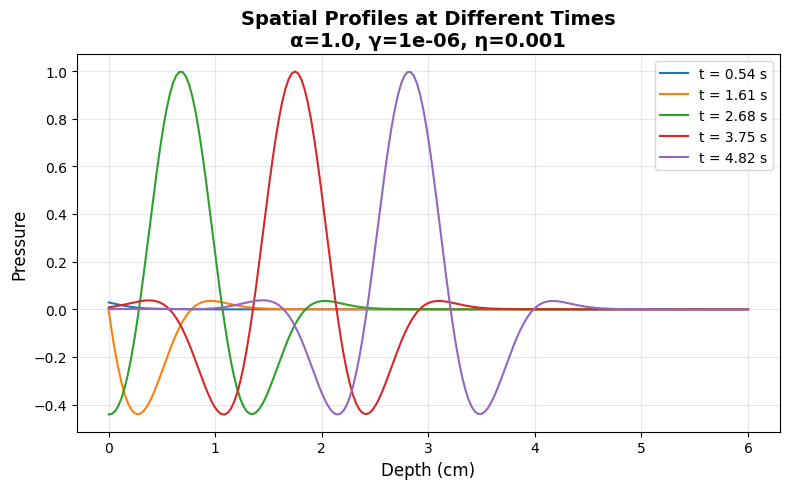}
		\caption{Waveform at Different Times}
		\label{fig:case1_time}
	\end{subfigure}
	\caption{Case 1: Waveform evolution in space and time.}
	\label{fig:case1_waveforms}
\end{figure}

Figures~\ref{fig:case1_loss} and \ref{fig:case1_propagation} present the loss component plots and the acoustic propagation field $p(x,t)$ for this case. Figures~\ref{fig:case1_depth} and \ref{fig:case1_time} reveal that the wave travels with very little distortion. The waveform retains its Gaussian‑modulated cosine shape across the entire domain. The waveforms are identical copy of the source pulse, just shifted later in time. The peak amplitude stays locked at $1.0$ across depth, and the trough stays at $0.43$ also. This means at low $\eta$, propagation is essentially non-distorting.

\begin{figure}[h!]
	\centering
	\begin{subfigure}[t]{0.45\textwidth}
		\centering
		\includegraphics[height=1.65in]{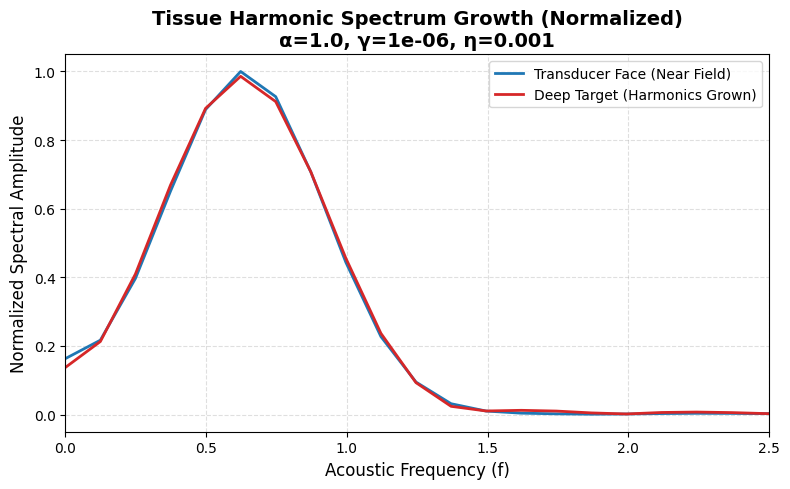}
		\caption{Tissue Harmonic Spectrum}
		\label{fig:case1_spectrum}
	\end{subfigure}%
	\begin{subfigure}[t]{0.45\textwidth}
		\centering
		\includegraphics[height=1.65in]{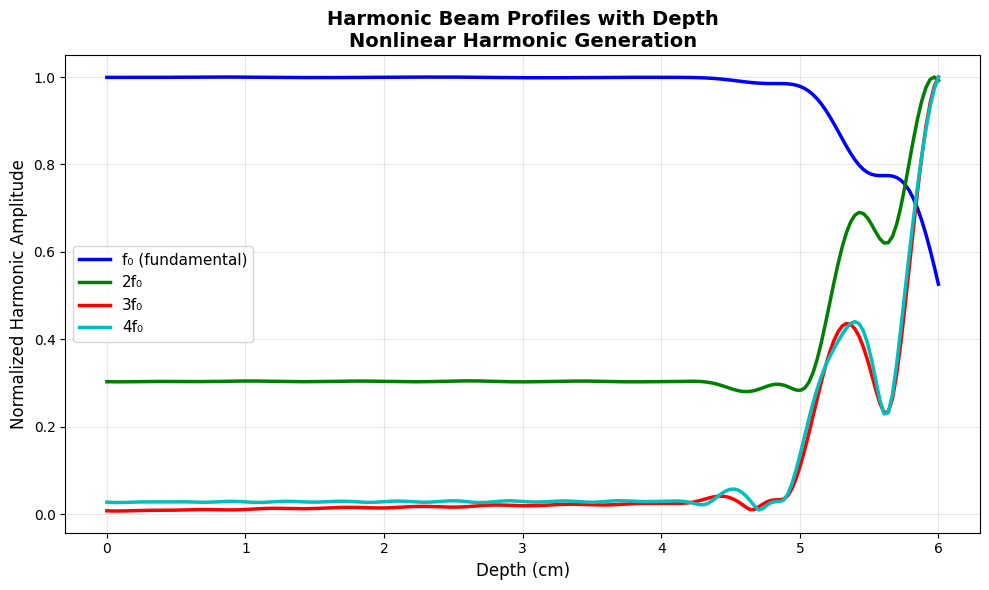}
		\caption{Harmonic Beam Profile}
		\label{fig:case1_beam}
	\end{subfigure}
	\caption{Case 1: Frequency analysis and beam profiles.}
	\label{fig:case1_freq}
\end{figure}
Figure~\ref{fig:case1_spectrum} shows that the near-field and deep-target curves are nearly the same across the entire frequency range. Both curves peak at about $f=0.65$, reaching a normalized amplitude of $ 1.0$ (near field) and about $0.99$ (deep target). Above $f \approx 1.3$, both curves collapse rapidly to near zero and stay flat out to $f=2.5$. This reveals that with low $\gamma$ and weak $\eta$, dissipative/diffusive losses are minimal. Meaning that no new high-frequency harmonic energy is generated or retained during propagation to depth. The deep target spectrum is just a damped copy of the near-field spectrum.

\begin{figure}[H]
	\centering
	\includegraphics[width=0.75\textwidth]{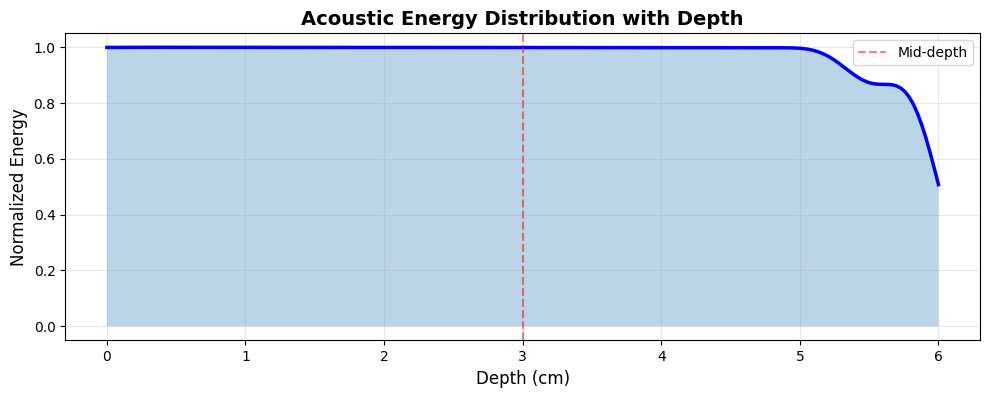}
	\caption{Case 1: Acoustic energy distribution}
	\label{fig:case1_energy}
\end{figure}

In Figure~\ref{fig:case1_beam}, the fundamental harmonic is flat at amplitude $ 1.0$ from $0 - 4.2$ cm, then a decline and rebound before a final drop to $ 0.52$ at $6$ cm. The second harmonic amplitude is stable at $0.30$ through the whole near field, dips slightly to $ 0.285$ around $4.6$ cm, then rises steeply starting at $4.9$ cm with a brief dip at $5.6$ cm before finally climbs to $1.0$ at $6$ cm. Both the third and fourth harmonic amplitude stay low and nearly flat $(0.01 - 0.03)$ through $4$ cm, with a small, shared bump around $4.5$ cm and a brief dip near $4.7$ cm. Then, both rise together after $4.9$ cm, peaking almost identically at $ 0.44$ around $5.3 - 5.4$ cm, dipping to $ 0.23$ at $5.6$ cm, and rise together to $1.0$ at $6$ cm. The acoustic energy distribution for this case is shown in Figure~\ref{fig:case1_energy}.

\subsection{Case 2: Stronger Nonlinearity, Low Attenuation}
\begin{figure}[h!]
	\centering
	\begin{subfigure}[t]{0.45\textwidth}
		\centering
		\includegraphics[height=1.65in]{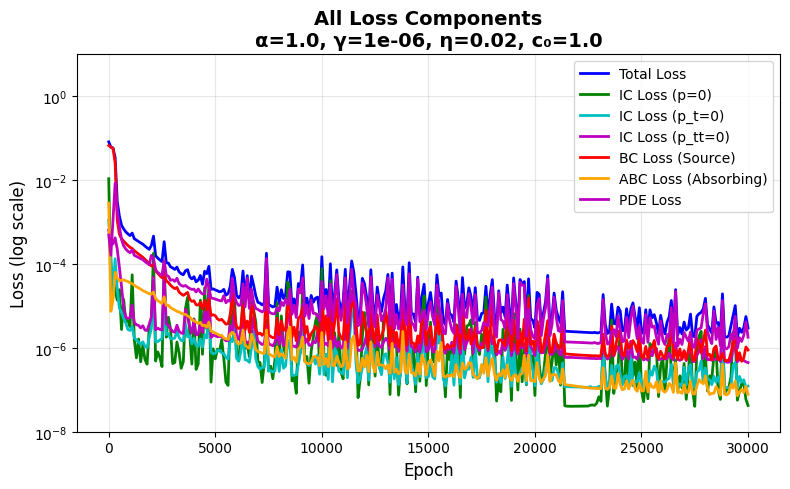}
		\caption{Loss Components}
		\label{fig:case2_loss}
	\end{subfigure}%
	\begin{subfigure}[t]{0.45\textwidth}
		\centering
		\includegraphics[height=1.65in]{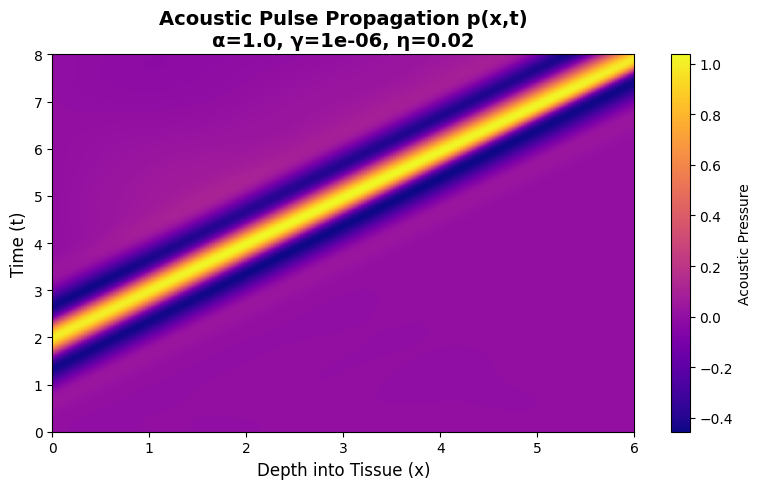}
		\caption{Acoustic Propagation}
		\label{fig:case2_propagation}
	\end{subfigure}
	\caption{Case 2: Loss components and acoustic propagation field.}
	\label{fig:case2}
\end{figure}

\begin{figure}[h!]
	\centering
	\begin{subfigure}[t]{0.45\textwidth}
		\centering
		\includegraphics[height=1.65in]{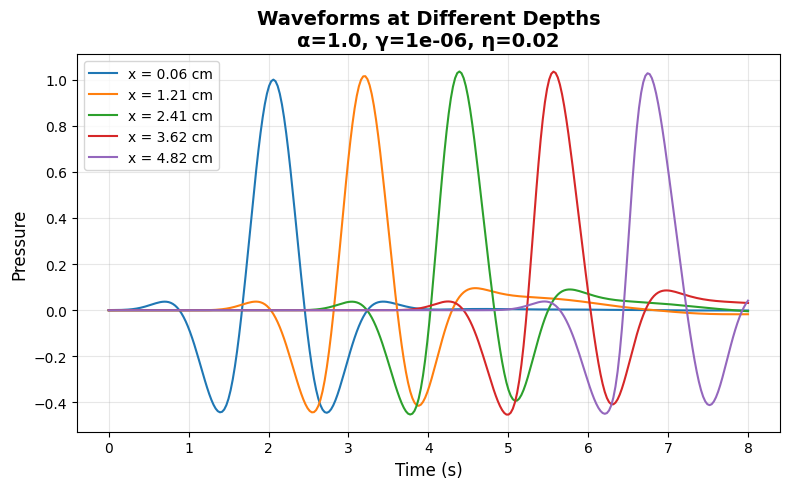}
		\caption{Waveform at Different Depths}
		\label{fig:case2_depth}
	\end{subfigure}%
	\begin{subfigure}[t]{0.45\textwidth}
		\centering
		\includegraphics[height=1.65in]{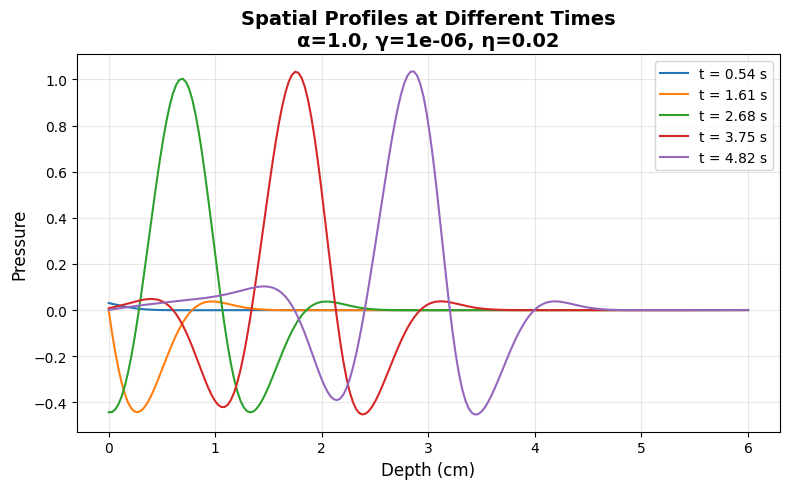}
		\caption{Waveform at Different Times}
		\label{fig:case2_time}
	\end{subfigure}
	\caption{Case 2: Waveform evolution in space and time.}
	\label{fig:case2_waveforms}
\end{figure}

Figures~\ref{fig:case2_loss} and \ref{fig:case2_propagation} present the loss component plots and the acoustic propagation field $p(x,t)$ for this case. In Figures~\ref{fig:case2_depth}, we observe from the source pulse (peak $1.0$, trough  $0.43$) that the waveforms change progressively with depth. The peak amplitude creeps above $1.0$, reaching roughly $1.03 - 1.04$ at the middle depths ($x=2.41 - 3.62$ cm) before reducing gradually. We notice also that a slowly decaying residual tail appears after each trough and grows somehow larger and longer-lived with depth.  This is due to the high value of $\eta$. Higher $\eta$ introduces cumulative, depth-dependent distortion. This reflects that energy is being fed into the waveform as it propagates (consistent with Figure~\ref{fig:case2_spectrum} below).  Figure~\ref{fig:case2_time} shows that waveform at different times.

\begin{figure}[h!]
	\centering
	\begin{subfigure}[t]{0.45\textwidth}
		\centering
		\includegraphics[height=1.65in]{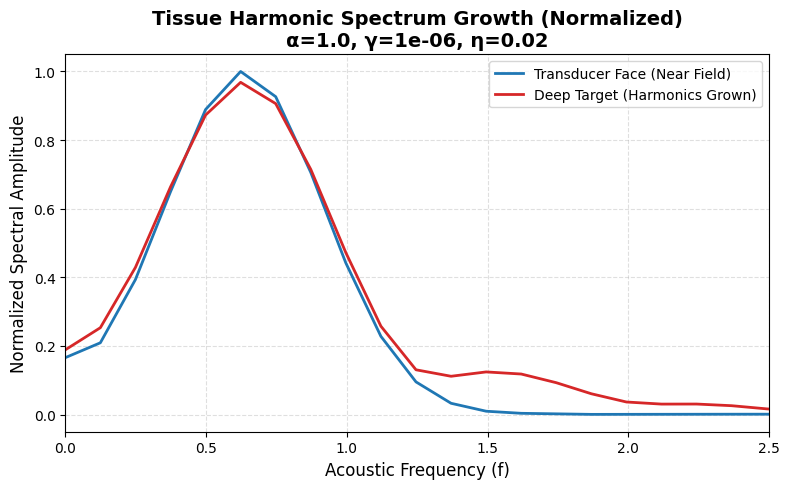}
		\caption{Tissue Harmonic Spectrum}
		\label{fig:case2_spectrum}
	\end{subfigure}%
	\begin{subfigure}[t]{0.45\textwidth}
		\centering
		\includegraphics[height=1.65in]{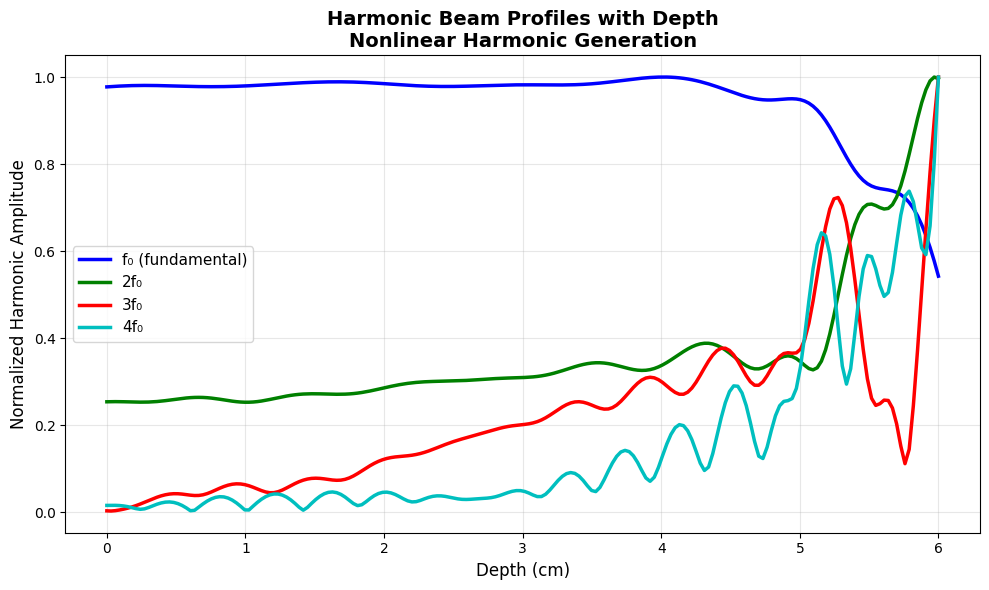}
		\caption{Harmonic Beam Profile}
		\label{fig:case2_beam}
	\end{subfigure}
	\caption{Case 2: Frequency analysis and beam profiles.}
	\label{fig:case2_freq}
\end{figure}
In Figure~\ref{fig:case2_spectrum}, the low-frequency behavior ($f < 1.1$) shows some minimal distortion unlike in Case 1. However, at higher frequencies $f > 1.1$, the distortion is notable. The deep-target curve (red) separates from the near-field curve (blue) and decays slower. The deep-target curve flattens into a broad plateau around $0.12–0.14$ spanning roughly $f=1.5$ to $f=2.5$, while the near-field curve continues down to $\approx 0.01–0.02$. This shows that increasing $\eta$ (with $\gamma$ still the same) does not change the fundamental/near-field spectrum much, but it substantially boosts and sustains higher-harmonic content at the deep target.

\begin{figure}[h!]
	\centering
	\includegraphics[width=0.75\textwidth]{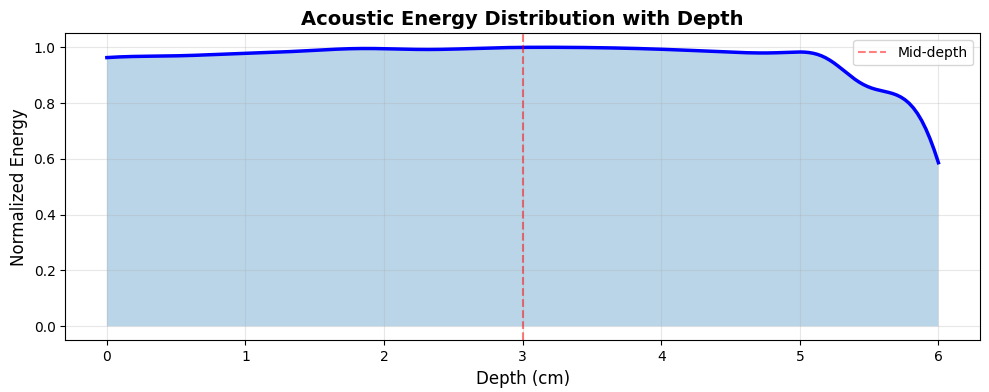}
	\caption{Case 2: Acoustic energy distribution}
	\label{fig:case2_energy}
\end{figure}

For this case, the fundamental harmonic, in Figure~\ref{fig:case2_beam}, starts around $0.98$, and rises to a peak of $1.0$ around $4$ cm. It then dips to $0.95$, and further declines with partial rebounds till it falls to $0.54$ at $6$ cm. The second harmonic is much noisier than in Case 1.  It oscillates between $0.25 - 0.31$ in the near field before trending upward, reaching $0.39$ near $4.3$ cm. It then dips and climbs again twice, before rising to $1.0$ at $6$ cm. The third and fourth harmonics show growing multiple oscillations with the fourth harmonic having the most multi-peak oscillatory behavior before they finally rise to $1.0$ at $6$cm. The acoustic energy distribution for this is presented in Figure~\ref{fig:case2_energy}.

\subsection{Case 3: Strong Nonlinearity, High Attenuation}
\begin{figure}[h!]
	\centering
	\begin{subfigure}[t]{0.45\textwidth}
		\centering
		\includegraphics[height=1.65in]{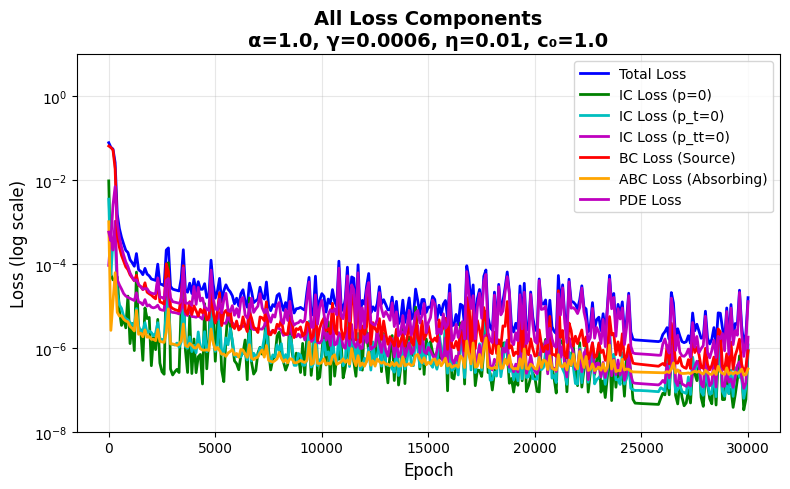}
		\caption{Loss Components}
		\label{fig:case3_loss}
	\end{subfigure}%
	\begin{subfigure}[t]{0.45\textwidth}
		\centering
		\includegraphics[height=1.65in]{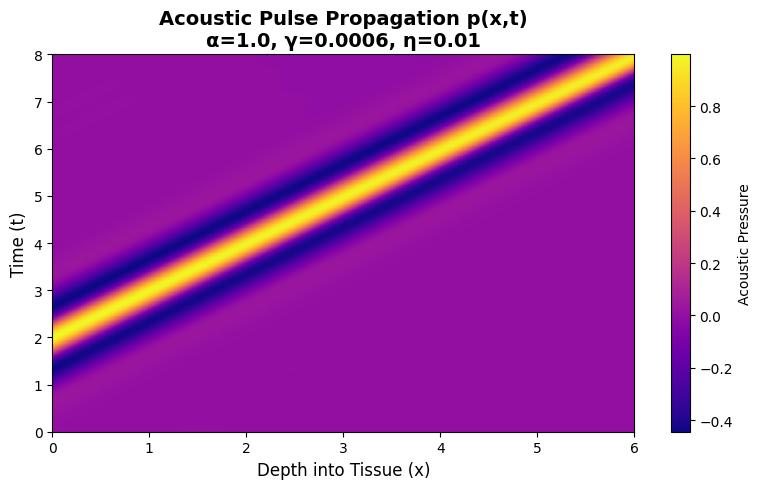}
		\caption{Acoustic Propagation}
		\label{fig:case3_propagation}
	\end{subfigure}
	\caption{Case 3: Loss components and acoustic propagation field.}
	\label{fig:case3}
\end{figure}

\begin{figure}[h!]
	\centering
	\begin{subfigure}[t]{0.45\textwidth}
		\centering
		\includegraphics[height=1.65in]{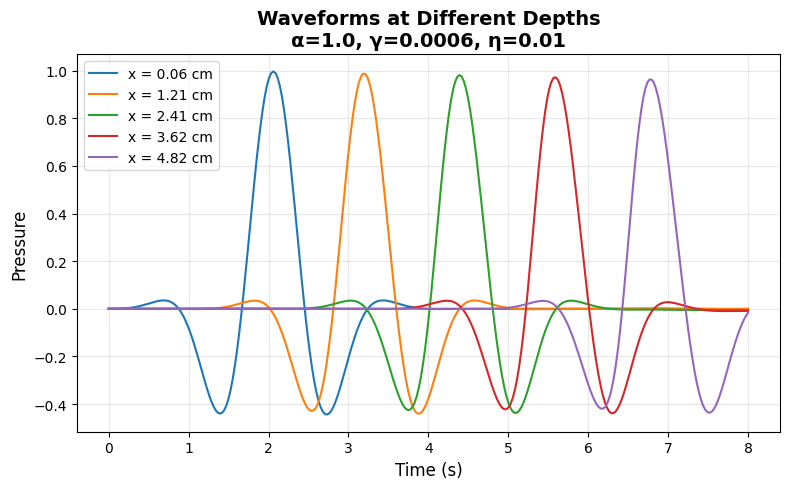}
		\caption{Waveform at Different Depths}
		\label{fig:case3_depth}
	\end{subfigure}%
	\begin{subfigure}[t]{0.45\textwidth}
		\centering
		\includegraphics[height=1.65in]{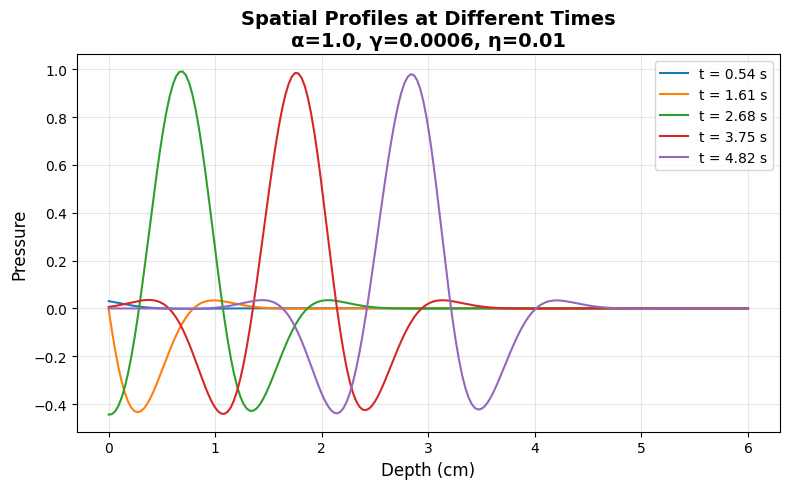}
		\caption{Waveform at Different Times}
		\label{fig:case3_time}
	\end{subfigure}
	\caption{Case 3: Waveform evolution in space and time.}
	\label{fig:case3_waveforms}
\end{figure}

Figures~\ref{fig:case3_loss} and \ref{fig:case3_propagation} present the loss component plots and the acoustic propagation field $p(x,t)$ for it.
This case sits between Case 1 and Case 2. In Figure~\ref{fig:case3_depth}, the source pulse, unlike Case 2, drifts gradually downward with depth. The trailing residual tail is present but much more modest than in Case 2. And it is not the same as in Case 1. This intermediate case behaves like a blend of cases 1 and 2.  It is observed that while $\eta$ drives distortion, $\gamma$ drives amplitude decay. With a less strong $\eta$ (0.01, half of Case 2) and high $\gamma$, the waveform is less distorted and the amplitude decays gradually. Figure~\ref{fig:case3_time} present the waveform propagation at different times.

\begin{figure}[h!]
	\centering
	\begin{subfigure}[t]{0.45\textwidth}
		\centering
		\includegraphics[height=1.65in]{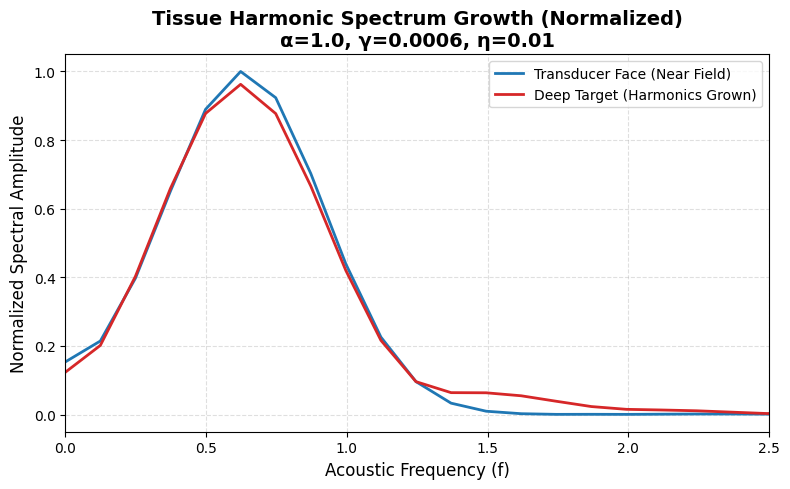}
		\caption{Tissue Harmonic Spectrum}
		\label{fig:case3_spectrum}
	\end{subfigure}%
	\begin{subfigure}[t]{0.45\textwidth}
		\centering
		\includegraphics[height=1.65in]{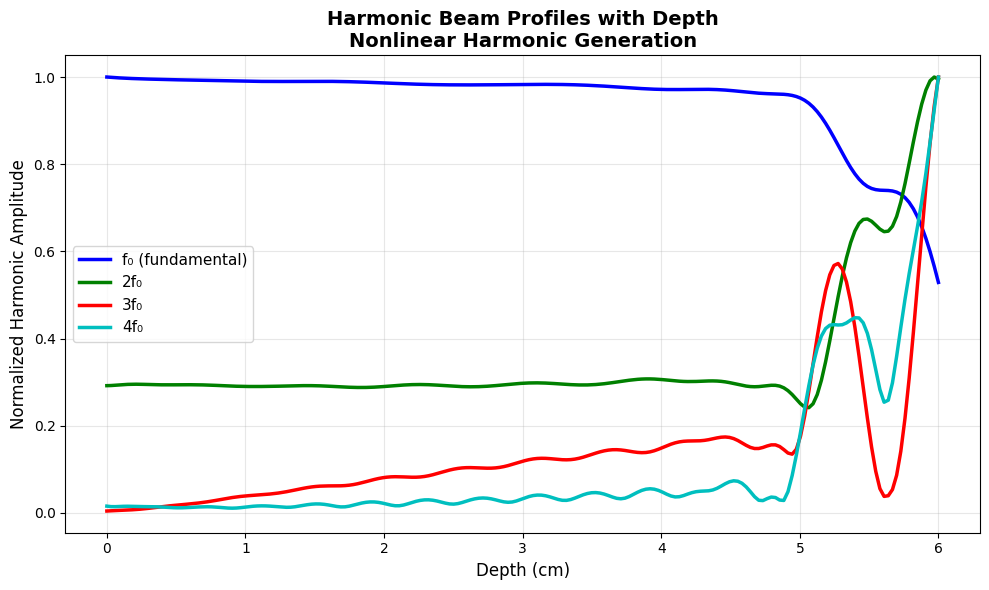}
		\caption{Harmonic Beam Profile}
		\label{fig:case3_beam}
	\end{subfigure}
	\caption{Case 3: Frequency analysis and beam profiles.}
	\label{fig:case3_freq}
\end{figure}

In Figure~\ref{fig:case3_spectrum}, the low-frequency onset is similar to case 1 but with a slightly wider gap. The peak is with amplitude of $ 1.0$ for near field and $0.93$ for deep target. And like Case 2, the deep-target curve pulls away from the near-field curve around $f\approx 1.21$, forming an elevated tail. However, the plateau here sits a bit lower, around $0.09 - 0.10$, as compared to Case 2. The result shows that with $\gamma$ being raised but $\eta$ reduced slightly as compared to Case 2, the elevated high-frequency tail at depth still appears. This reveals that this behavior is governed primarily by $\eta$ rather than $\gamma$.  However, the magnitude is smaller due to the increase in the value of $\gamma$. 

\begin{figure}[H]
	\centering
	\includegraphics[width=0.75\textwidth]{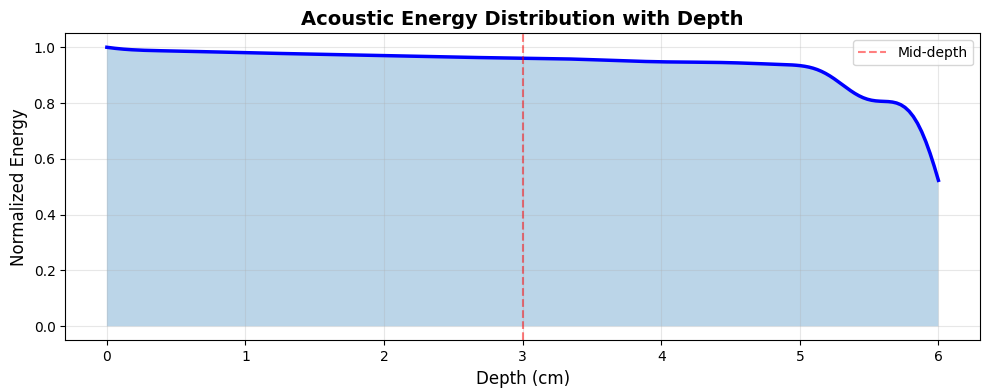}
	\caption{Case 3: Acoustic energy distribution}
	\label{fig:case3_energy}
\end{figure}

Here, the fundamental harmonic (see Figure~\ref{fig:case3_beam} ) starts at $1.0$, declines gradually and smoothly (less oscillation than Case 2) to $0.95$ by $5$ cm, then drops with partial rebound before a final drop to $0.52$ at $6$ cm. The second harmonic starts about $0.20$ only mild oscillation. It dips to $0.24$ near $5$ cm. It then rises and dips again before finally climbs to $1.0$ at $6$ cm. The third harmonic amplitude rises more steadily and smoothly than in Case 2. It reaches $0.17$ by $4.5$ cm and dips to $0.13$ near $4.9$ cm before surging to a peak of $0.57$ at $5.3$ cm. It crashes to $0.04$ at $5.75$ cm (the deepest dip of the three cases) and a final rise to $1.0$. The fourth harmonic amplitude has mild oscillation through the near field (staying under $0.07$ until $4.5$ cm).  It experiences a small dip near $4.8$ cm and another dip at 5.6 cm before a steep rise to $1.0$ at $6$ cm. The acoustic energy distribution for this case is shown in Figure~\ref{fig:case3_energy}.

\section{Conclusion}
\label{conclude}
We solved the Westervelt equation \eqref{eq:bc_full_absorbing} with PINNs for three synthetic cases to separate regimes that overlap in real tissues (Water/Brain/Fat $\eta=0.0035-0.008$ within $2\%$ at 6 cm depth). PINNs enables direct enforcement of the nonlinear wave physics. It requires neither mesh generation nor traditional numerical discretization. Through the study of three cases, we demonstrate the PINNs' capability to capture some nonlinear phenomena. These include harmonic generation and waveform steepening.

The near-field spectrum is invariant across all three parameter sets. $\eta$ controls whether harmonics survive/grow at depth. At $\eta=0.001$, deep-target and near-field curves are nearly identical (no harmonic growth). At $\eta=0.02$ and $\eta=0.01$, a broad elevated curve appears in the deep-target spectrum above $f \approx 1.1 - 1.5$ (as seen in cases 2 and 3), indicating $\eta$ drives sustained harmonic generation/energy transfer to higher frequencies as the wave travels to the deep target. However, the value of $\gamma$ appears to add an effect in lowering the magnitude of the rise as shown in Case 3.

The study has some limitations. These suggest directions for future research. First, we use one‑dimensional geometry, which is good for fundamental understanding. However, it does not capture the full three‑dimensional complexity of clinical HIFU applications. Extending the study to two or three dimensions would enable simulation of these applications. Second, we normalized the parameters for demonstrating the method. Yet they do not correspond directly to physical tissue properties. Future work will incorporate tissue parameters based on experimental measurements. Third, we assume homogeneous media. Incorporating spatially varying parameters would enable simulation of tissue heterogeneity especially for clinical applications. 

The Loss components reached at least $10^{-5}$ across all cases. These results demonstrate that $\gamma$  scales the rate of energy absorption (amplitude decay), while $\eta$ governs frequency-dependent phase speeds (dispersion). These two parameters dictate whether a nonlinear pressure pulse can successfully transfer energy through deep fluid layers or if it degrades into low-frequency environmental noise. This is useful for both diagnostic and therapeutic ultrasound applications.

%\newpage
\bibliographystyle{unsrt}
\bibliography{refWestervelt}

\end{document}